\documentclass[12pt]{article}
\usepackage[hmargin=1in, vmargin=1in]{geometry}

\usepackage{amsmath}
\usepackage{amsthm}
\usepackage{amssymb}
\usepackage{hyperref}
\usepackage{caption}
\usepackage{float}
\usepackage{color}

\usepackage{graphicx}
\usepackage{dsfont}
\usepackage{indentfirst}
\usepackage{sidecap}
\usepackage{multicol}
\usepackage{gensymb}
\usepackage{subcaption}

\begin{document}




\title{Talk to the Hand:  Generating a 3D Print from Photographs}

\author{Edward Aboufadel, Sylvanna V. Krawczyk, \\ and Melissa Sherman-Bennett\footnote{EFA:  Grand Valley State University; SVK:  California State University (Sacramento), MS: Bard College at Simon's Rock.  Contact: aboufade@gvsu.edu }}

\maketitle



\abstract{This manuscript presents a linear algebra-based technique that only requires two unique photographs from a digital camera to mathematically construct a 3D surface representation which can then be 3D printed. Basic computer vision theory and manufacturing principles are also briefly discussed.}






\section{Introduction}

\subsection{Purpose}

\begin{figure}[h]
        \centering
        \begin{subfigure}[b]{0.25\textwidth}
                \includegraphics[width=\textwidth]{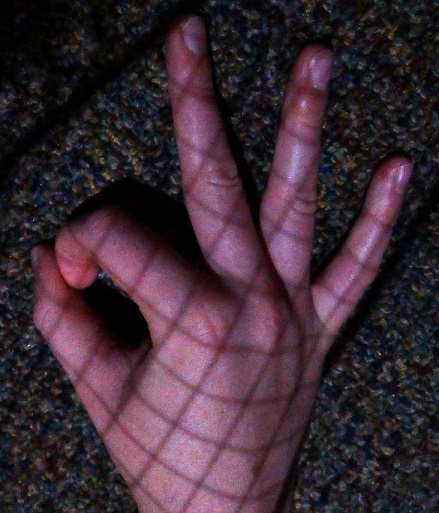}
        \end{subfigure}
		\qquad
        \begin{subfigure}[b]{0.25\textwidth}
                \includegraphics[width=\textwidth]{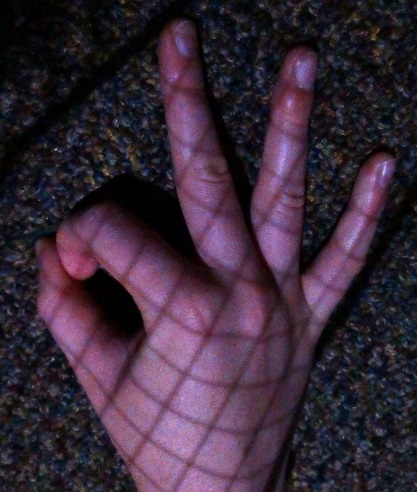}
        \end{subfigure}
		\qquad
        \begin{subfigure}[b]{0.25\textwidth}
               \includegraphics[width=\textwidth]{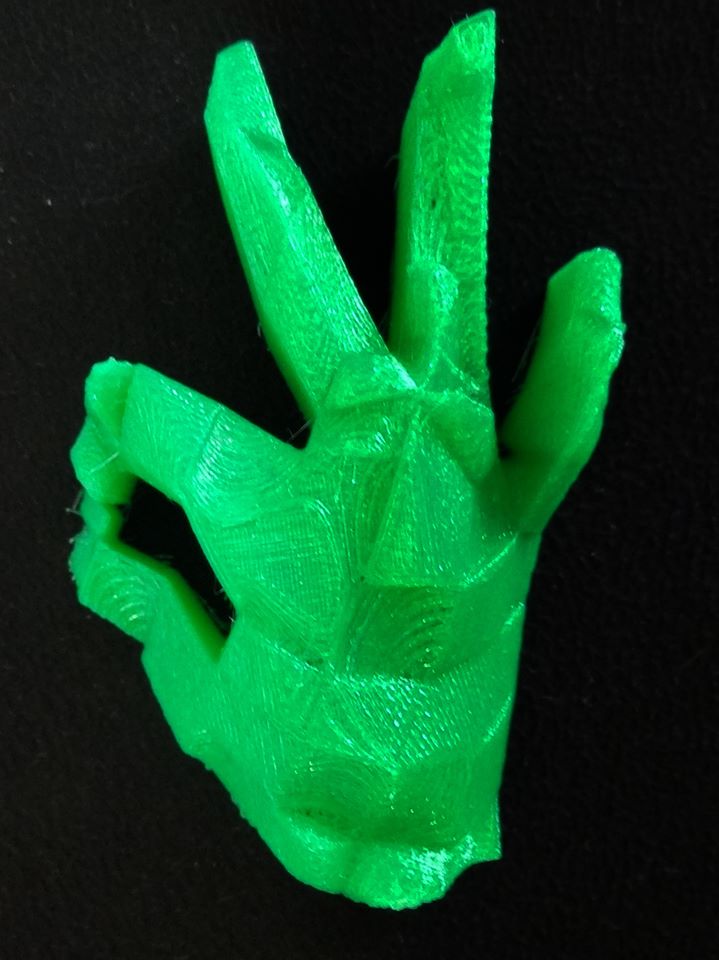}
        \end{subfigure}
        \caption{3D print (c) generated from two photographs (a, b).}
        \label{fig:handpic}
\end{figure}

The purpose of this paper is to describe an algorithm based on ideas from linear algebra and geometry to generate a 3D-printed model from a pair of similar photographs (see Figure \ref{fig:handpic}). Our work is part of a broader movement to answer the question, ``What can a mathematician do with a 3D printer?'' for which there is some recent work done in areas such as fractals \cite{Hart}, knots \cite{MakerHome}, smooth manifolds \cite{Segerman}, polyhedra \cite{MAA, MakerHome}, a print based on Conway's Game of Life \cite{Conway}, and demonstrations of theoretic or historical constructs \cite{Knill}.  The recent availability of relatively inexpensive 3D printers has facilitated these projects.

3D printers create objects through the technology of additive manufacturing -- starting with nothing and building up the object layer by layer, using only the material needed for the object itself.  The most affordable 3D printer models are extrusion-based models, where a thin plastic filament is fed into a heating element called an \emph{extruder} which melts the filament and lays down a thin trail of plastic onto the build plate.

The opportunities these printers present for at-home manufacturing naturally lead to the problem of printing replicas of an already-existing object. A number of websites and apps, such as \cite{Scanner} and \cite{AutoDesk}, utilize one approach to this problem, which involves taking many photographs of an object at varying angles. A high-resolution, 360\degree, 3D model of the object is produced from these photographs, which can then be printed. The drawbacks of this approach is that it requires a large number of photographs, often over 20, and does not function well on objects with smooth, homogeneous surfaces. We present an alternative method to 3D print models of existing objects. It is also photography-based, but utilizes only two photographs, has minimal restrictions on the nature of the object to be modeled, and is simple enough to be implemented using mathematical software packages such as \emph{Mathematica}. The method produces a bas-relief model, which is sufficient for our purposes.

We will begin with more information about 3D printing in general.  Then, we will discuss the theoretical basis of creating 3D models from photographs via the recovery of depth data. Finally, we will describe specifically how the hand in Figure \ref{fig:handpic} was generated from indicated photographs.

\subsection{About 3D Printing}

Inexpensive 3D printers that are now available typically create objects made of polylactide (PLA) or acrylonitrile butadiene styrene (ABS) plastic filament.  The quick-to-solidify melted plastic is extruded onto a flat surface called the ``build plate''.  The extruder moves in the \emph{xy}-plane parallel to the build plate. Once the extruder has deposited one layer of plastic onto the build plate, the build plate moves down slightly along the \emph{z}-axis and the extruder lays a new horizontal layer on top of the previously extruded material. In this way, the printer stacks hundreds of horizontal layers on top of each other to manufacture user-designed objects.

Computers represent 3D objects using a mesh, which is a collection of polygons in $\mathds{R}^3$ that share edges, forming one continuous surface. Both the type and number of the polygons on the surface can be altered to change the resolution of the object. The format of files containing information about meshes varies depending on the type of polygon used. A common file format is STereoLithographic (STL), which handles triangle meshes. STL files record the vertices and normal vectors of each triangle, and can be either binary (computer-readable) or ASCII (human-readable). For examples, see \cite{MAA} and \cite{Primer}. STL files, and other polygon mesh files, are transformed into instructions to the printer on how to print an object (see \cite{Primer}), so once an STL representation of an object has been created, the object can be printed. STL files can be written using a text editor, or generated from a 3D designer program.  In addition, \emph{Mathematica} can create an STL file for any 3D surface it can display.  STL files representing a mathematically defined surface, our particular area of interest, can be created by defining a function $z=f(x,y)$ to represent the height $z$ of a surface at points $(x,y)$ in a bounded region and then graphing this function and exporting the mesh of its graph as an STL file. For example, altitude data (as a function of longitude and latitude) can be used to create a mesh, and STL file, for a mountain. The function will not necessarily be differentiable, but in order to print properly, we want it to be continuous.

Some recent papers connecting mathematics and 3D printing include \cite{MAA}, \cite{Primer}, \cite{Austin14}, \cite{Segerman}, and \cite{Escher}.


\section{Experimental Section:  Underlying Theory On Extracting Depth Data}

There are two main techniques for finding the 3D shape of an object from photographs: \emph{photometric} and \emph{geometric}. Photometric methods use the manner in which different parts of the object reflect light to derive the orientation of parts of the object's surface, by applying different ``shading functions"  \cite{springer}. These functions describe how light will bend when it hits a specific material. However, to use the shading functions to recover the shape of an object, one must know exactly the type and direction of light used in the photo. As this is impractical for most photographs taken outside of a photography studio, we instead utilized a geometric approach, which relies upon the nature of photography as a projection from $\mathds{R}^3$ to $\mathds{R}^2$, or from the real world to the ``image plane'' of the flat photo. If one considers more than one photograph, this projection produces a number of similar triangles, from which depth can be recovered in a process aptly called \emph{triangulation}.

\subsection{Perspective Projections}


An interesting article on determining camera location from a photograph, which is connected to our work, is \cite{Camera}. As described in detail in that article, the projection that occurs when an object is photographed is not linear, that is, it is not as simple as mapping $(x, y, z)$ to $(x, y, 0)$.  Instead, in what is called a \emph{perspective transformation}, a given point is projected onto the image plane along the ray connecting the point to the camera lens. ``Taking a picture'' maps a point $(x, y, z)$ in 3-space to its perspective projection $P$ by the mapping:

$$ P: (x, y, z) \mapsto (-f\frac{x}{z}, -f\frac{y}{z}, 0),$$  where $f$ is the focal length, or the distance from the camera lens to the image plane.  (See also \cite{springer}.)  Under this projection, proportion and angle are preserved only for lines and shapes in a plane parallel to the camera (and thus to the image plane). All others are warped, with points further from the camera being affected less.

This mapping provides one way to recover the depth of points: if the original $x$- and $y$- coordinates of a point in $\mathds{R}^3$, the focal length of the camera, and the point's location in the image plane (which we will call its ``pixel coordinates'') are all known, the depth $z$ can be found easily. However, the significant external measurements required render this technique as impractical as the photometric methods. We will instead utilize a convenient relationship that arises in perspective projections from slightly different positions, which relates depth to a quantity easily measurable from the photos themselves.

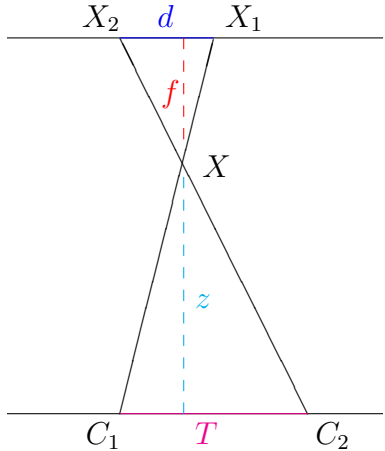
\begin{figure}

\centering
\setlength{\unitlength}{5cm} 
\begin{picture}(1, 1.1)  

\put(0, 0.1){\line(1,0){1}}
\put(0, 1.1){\line(1,0){1}}
\put(0.3, 0.1){\line(1,4){0.249}}
\put(0.8, 0.1){\line(-1,2){0.499}}
\put(0.21, 0.02){$C_1$}
\put(0.82, 0.02){$C_2$}
\put(0.2, 1.13){$X_2$}
\put(0.58, 1.13){$X_1$}
\put(0.52, 0.73){$X$}

\color{blue}
\put(0.4, 1.13){$d$}
\put(0.3, 1.1){\line(1,0){0.25}}

\color{magenta}
\put(0.5, 0.02){$T$}
\put(0.3, 0.1){\line(1,0){0.5}}

\color{cyan}
\put(0.5, 0.38){$z$}
\multiput(0.47,0.1)(0,0.06){11}
{\line(0,1){0.03}}

\color{red}
\put(0.41, 0.93){$f$}
\multiput(0.47,1.1)(0,-0.06){5}
{\line(0,-1){0.03}}

\end{picture}

\caption{$\Delta X_1X_2X \sim \Delta C_1C_2X$, therefore $z=\frac{fT}{d}$, where $z$ is depth, $f$ is the focal length, $X_1$ is the projection of $X$ by Camera 1 ($C_1$) and $X_2$ is the projection of $X$ by Camera 2 ($C_2$).  The translation between the cameras is $T$, while the translation between the projection of $X$ is $d$. Image derived from \cite{cheatsheet}.}
\label{fig:simTriangles}
\end{figure}

\subsection{Depth-Disparity Relationship}

We consider two photographs with parallel image planes, meaning the position of the camera in the two photos differs only by a translation $T$. For our algorithm, the translation $T$ must be parallel to the image plane, not perpendicular to it. Figure \ref{fig:simTriangles} is a simplified sketch of such a situation. (A more detailed image and discussion can be found in \cite{cheatsheet}.) One can see that in such a situation, $d$, or the distance between a point's projection in one image plane and its projection into the other, depends on how close the point is to the camera. In fact, the closer the point is to the camera, the larger $d$, which we will call the ``disparity," will be. So, we observe that:
$$z \propto \frac{1}{d}$$
where $z$ is the depth of a point. Geometric reasoning with similar triangles or algebraic manipulations of the locations of the point on the two image planes confirms this relationship, as well as providing the constant of proportionality. Since 3D modeling only requires relative depth of points, and the constant of proportionality involves quantities we may not know, the above proportion is sufficient for our purposes.


\section{Results and Discussion: Creating and Printing the Hand}

Given the preliminary ideas above, we will now describe how to extract depth data from two photographs. After a general discussion of our algorithm, we will explain how it is used with a specific example. For simplicity's sake, we will describe how this is done with a contrived scene, using two pictures of a cube with side length of 1 inch (see Figure \ref{fig:tmapstos}).  Then, we will show how to apply our method to create and print a hand, and suggest directions for future work.

\subsection{The General Idea}

Our method, implemented in \emph{Mathematica}, is a geometric approach that uses two photographs of an object, taken from camera positions so that the translation of the camera, $T$ in Figure \ref{fig:simTriangles}, is relatively small. It is based on the inverse relation between depth and disparity described in the section above; namely, given two photographs with image planes differing only by a small translation, one can recover the depth of a 3-D point by comparing where it is projected, or its \emph{pixel coordinates}, in the two image planes.

The key idea is that we only extract coordinates from the digital images for a pre-defined, relatively small set of points: the corners of approximately planar quadrilateral sections of the surface. We then use these coordinates and linear algebra techniques to compute the depths of the points within these quadrilateral faces. The calculated depths of individual quadrilateral-shaped surfaces are then combined to yield the complete depth data of the surface in question, and to create a depth map of the surface. Once this is done, the depth map can be converted easily into an STL mesh, which can then be printed.

\subsection{The Transformation Matrix} \label{findtransform}

\begin{figure}
\begin{center}
\includegraphics[width=6in]{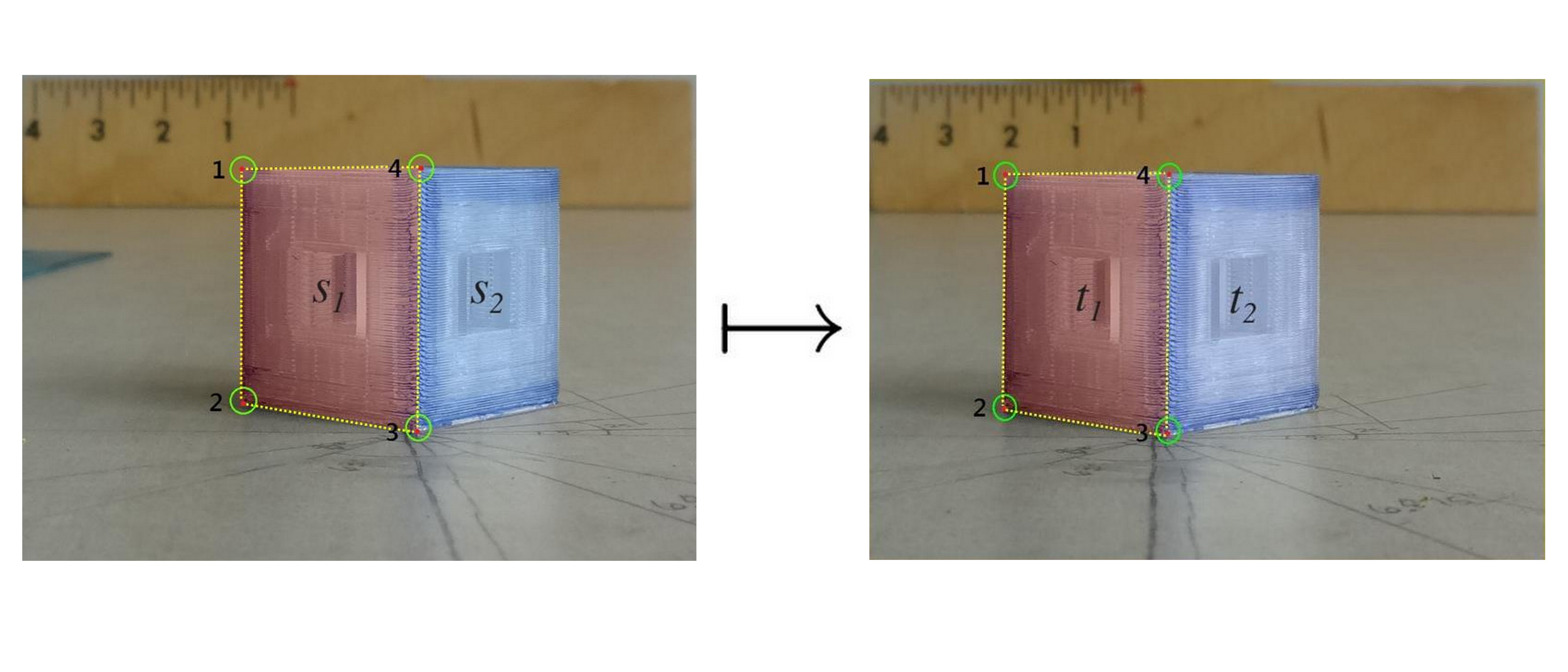}
\caption{Two example photographs, with corresponding planes \textbf{s} and \textbf{t} highlighted. The four pairs of corresponding pixel coordinates are also marked.}
\label{fig:tmapstos}
\end{center}
\end{figure}

For a plane $\textbf{s}$ in one image and its corresponding plane $\textbf{t}$ in the other image (see Figure \ref{fig:tmapstos} for an example), one can find \cite{Austin13} a transformation matrix $\textbf{M}_{\textbf{st}}$ and a constant $w$ such that, for any point with pixel coordinates $(x, y)\in\textbf{s}$ and pixel coordinates $(u, v)\in\textbf{t}$,
\begin{align}
\begin{bmatrix}
x & y & 1
\end{bmatrix}
&\textbf{M}_{\textbf{st}}
=w
\begin{bmatrix}
u & v & 1
\end{bmatrix}\\
\begin{bmatrix}
x & y & 1
\end{bmatrix}
& \begin{pmatrix}
a & d & g\\
b & e & h\\
c & f & j
\end{pmatrix}
= w
\begin{bmatrix}
u & v & 1
\end{bmatrix} \label{transform}
\end{align} where $w=gx+ hy +j$.
The elements of $\textbf{M}_{\textbf{st}}$ can clearly be determined if $(x, y)$ and $(u, v)$ are known. For one pair of corresponding pixel coordinates $(x_i, y_i)$ and $(u_i, v_i)$, rearranging Eq. \ref{transform} yields a 2$\times$9 matrix equation with $\{a, ..., j\}$ as unknowns:
\begin{align}
\begin{bmatrix}
x_i & y_i & 1 & 0 & 0 & 0 & -x_i u_i & -y_i u_i & -u_i\\
0 & 0 & 0 & x_i & y_i & 1 & -x_i v_i & -y_i v_i & -v_i
\end{bmatrix}
\begin{bmatrix}
a\\
b\\
c\\
d\\
e\\
f\\
g\\
h\\
j
\end{bmatrix}
=
\begin{bmatrix}
0\\
0
\end{bmatrix} \label{coeff}
\end{align}

If four distinct pairs of corresponding pixel coordinates $(u_i, v_i)$ and $(x_i, y_i)$ are used, four different versions of Eqn.~\ref{coeff} are produced and can be combined into one 8$\times$9 matrix, the nullspace of which yields the coefficients $\{a, ..., j\}$, up to a multiplicative constant. The scaling constant $w_i = g x_i + h y_i + j$ is determined separately for each point $(x_i, y_i) \in \textbf{s}$ after the elements of $\textbf{M}_{\textbf{st}}$ are found. Because the system is homogeneous, it is consistent, and the implication of the discussion in \cite{Austin13} is that the rank of the 8$\times$9 matrix will be 8 as long as the four points are not collinear.  So, we will have only one free variable.

\subsection {The Cube}

A 3-D model created from two photographs of a cube (see Figure \ref{fig:tmapstos}) serves to illustrate the basic process of the algorithm; models of more complex objects require more computation in each step, but do not require additional steps.  The cube in the photograph must be partitioned into quadrilateral planar faces with corners distinguishable in both photographs, as shown in the figure. The pixel coordinates of the corners in one photograph must be recorded and matched with those of the corresponding corner in the other photograph.  In the current algorithm, the user chooses the quadrilateral planar faces and finds the pixel coordinates of the corners manually. The latter can be done using any program that displays the pixel coordinates of a given pixel, including \emph{Microsoft Paint} and \emph{Mathematica}.

\subsubsection{Creating the depth map for one quadrilateral face}

\begin{figure}
\begin{center}
\includegraphics[height=1.8in]{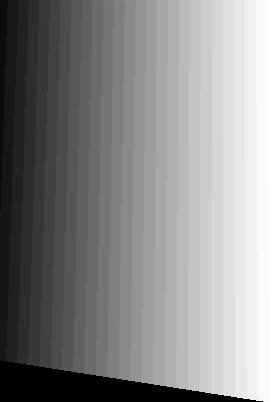}
\caption{A visual representation of the depth array of $\textbf{s}_1$, with values in the array interpreted as grayscale values (how light or dark a pixel is). Points closest to the camera are closer to white, while points further away are closer to black.}
\label{fig:leftmapgray}
\end{center}
\end{figure}

After the quadrilateral faces have been identified and the pixel coordinates of corresponding corners have been found, the transformation matrix $\textbf{M}_{\textbf{st}}$ for one face of the cube can be calculated as described in section \ref{findtransform}. We will here consider the left-hand face of the cube, choosing to find the matrix that maps $\textbf{s}_1$ to $\textbf{t}_1$, though the inverse operation is equivalent in this context. Thus, for any pixel coordinate $(x,y)\in \textbf{s}_1$, its corresponding pixel coordinate $(u,v) \in \textbf{t}_1$ is given by the equation

\begin{align}
\frac{1}{g\text{x} + h\text{y} + j}
\begin{bmatrix}
\text{x} & \text{y} & 1
\end{bmatrix}
\begin{pmatrix}
a & d & g\\
b & e & h\\
c & f & j
\end{pmatrix}
=
\begin{bmatrix}
\text{u} & \text{v} & 1
\end{bmatrix}\label{transform2}
\end{align}

The pixel coordinates $(x,y)$ and $(u,v)$ are the projections of a point \textbf{p} on the cube into the image planes of the respective photographs. The true depth $D$ of the point \textbf{p} is

\begin{align}
D = \frac{k}{\sqrt{(x-u)^2 + (y-v)^2}} \label{depth}
\end{align}
where the denominator of the fraction is the disparity and $k$ is an unknown constant of proportionality. Because the 3-D model will be a scaled version of the cube, however, relative depth is sufficient, and the value of $k$ can be chosen conveniently.
Finding the depth of each point on the left-hand face of the cube is simply a matter of calculating Equation \ref{depth} for each pixel coordinate in $\textbf{s}_1$. Note that this depth is measured from the camera, so the depths of corners 3 and 4 would be smaller than those of corners 1 and 2. Given the bottom-up manner in which our model will be printed, it is preferable to have the depth of points from the ``back" of the cube instead. To do this, we choose $k$ so that $D<1$ for all $(x, y) \in \textbf{s}_1$, and record $1- D$ rather than $D$.

\begin{figure}
\begin{center}
\includegraphics[width=2in]{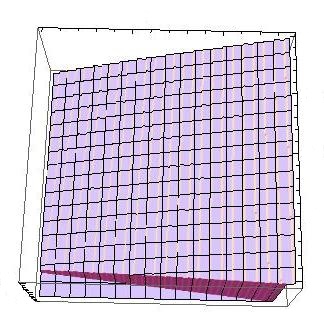} 
\caption{The 3-D version of the depth array in Figure \ref{fig:leftmapgray}, with values of the array interpreted as heights.}
\label{fig:leftmap3D}
\end{center}
\end{figure}

\emph{Mathematica} is capable of interpreting the depth array in a number of ways, two of which are shown in Figures \ref{fig:leftmapgray} and \ref{fig:leftmap3D}. In both, one can see that points outside of $\textbf{s}_1$ have a depth of zero, the line between corners 3 and 4 is closest to the camera, and the line between corners 1 and 2 is the furthest away, as expected. The vertical bars in which the depth remains the same, reflected either as constant grayscale value or constant height, are due to rounding errors. The slight jaggedness of the edges of the face, as well as the lines separating different areas of constant depth, is due to the discrete nature of the array. The 3-D version of the array, or the depth map, is what will eventually be used to create the mesh for our model of the cube. The grayscale version is used in the intermediate steps of the process for troubleshooting, as it is faster to generate and somewhat easier to manipulate.

\subsubsection{Creating the complete depth map}
\begin{figure}[h]
        \centering
        \begin{subfigure}[b]{0.4\textwidth}
                \includegraphics[width=\textwidth]{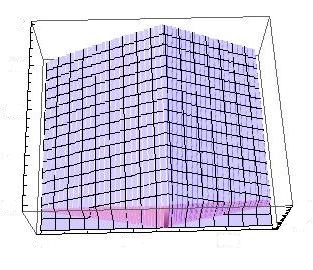}
        \end{subfigure}
		\qquad
        \begin{subfigure}[b]{0.4\textwidth}
                \includegraphics[width=\textwidth]{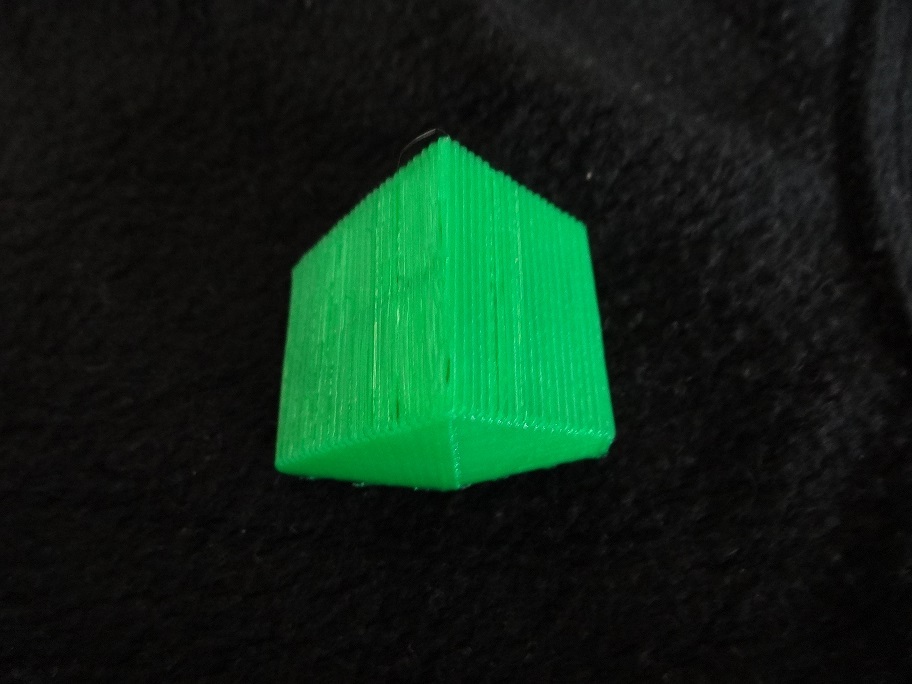}
        \end{subfigure}
        \caption{Left: The complete depth map for the cube, which is used to create the model.  Right: A 3D printed version of the depth map (a).}
\label{fig:depthmap3D}
\end{figure}

Creating the complete depth map for the cube is relatively simple. The process described above is duplicated for the right-hand face of the cube, resulting in two depth arrays. The constant $k$ must be the same for both arrays, or the relative depths of the two faces will be inaccurate. Once the two arrays are computed, they are combined to create one large, complete depth array. This step is mostly a matter of matching corners, as the upper-right corner of the left face is the same as the upper-right corner of the right face. Since the pixel coordinates of this corner in both photographs are known, it is possible to determine the position of that corner's depth values in each depth array. When combined, the depth arrays of the two photos must ``match up" at this position. Given this information and the dimensions of both depth arrays, each array is padded with columns or rows of zeros so that their shared corners are in the same position in each array, which ensures that the shared edge of the faces is indeed shared in the 3D model. The arrays are then simply added element-wise. After the depth arrays have been combined, a 3-D version of the large array is generated and exported from \emph{Mathematica} as a mesh, in the STL file format. It is then printed, creating the model shown in Figure \ref{fig:depthmap3D}. Note that the 3-D model of the cube is accurate to the cube as it appears in the photographs, not as it was originally. Distortions occurring due to the perspective projection involved in taking a photograph are not corrected.

\subsection {The Hand}

Moving from photographs of a cube to the photographs of the hand, the first alteration in the algorithm is imposed by the fact that hands, unlike cubes, are not comprised of planar faces with distinct edges and corners. As hands are relatively smooth and featureless, it was necessary to alter the hand during the process of taking the photograph in order to create distinct corners for the quadrilateral faces. As shown in Figure \ref{fig:handpic}, a shadow grid was projected onto the hand by shining a common 1" diameter LED flashlight onto a grid printed on a transparency. The matching pairs of pixel coordinates were then found by the user, as with the cube.

This additional step has the advantage of clearly delineating quadrilateral sections of the hand. However, the quadrilateral regions are not necessarily planar; in assuming that they are, as the rest of the algorithm demands, one necessarily produces an approximation of the hand as it appears in the photograph. Any curvature within the quadrilaterals will not be reflected in the model, as can be seen with the knuckle of the index finger. This can be mitigated by using a grid with more squares per inch, which has the effect of increasing the accuracy of the approximation. Additionally, since it only takes 3 points to determine a plane and the transformation matrix requires 4 coplanar points, a smaller grid decreases the likelihood that one of the corners of the quadrilaterals is not coplanar with the others.

\subsubsection{Creating the complete depth map}

The process of combining the individual depth maps into the complete hand depth map is more involved than that for the cube simply because there are more quadrilateral faces to consider. The method used to join the depth maps making up the cube can be used for the depth arrays of any two adjacent faces; if two depth arrays are combined, the resulting array can then be combined with another depth array, providing the two faces have an edge in common. Thus, the process of creating the hand depth map is that used to create the cube depth map, but applied iteratively. The requirement that depth arrays must share an edge to be combined imposes an order on the combination process. The ``path" of combination can be any path that results in all of the arrays being combined; if multiple such paths are available, they yield identical results.

\subsubsection{Correcting the depth map}

The larger number of faces introduces some error in the depth map, which must be corrected before printing the model. The boundaries of the quadrilateral faces are determined by lines, and the arrays, which are discrete, can only approximate these lines. For this reason, some pixels on the boundaries may be interpreted as being in multiple quadrilateral faces, and their depth values are stored in multiple arrays. When these redundancies occur, the arrays overlap slightly when they are added, resulting in pixels with comparatively large values. In the depth map, this produces large, very thin spikes, which interfere with the printing process. This is mitigated by replacing these depth values by the average of nearby values, or by manually altering the depth array.


\section{Conclusions and Future Directions}

The main limitation of the algorithm is the level of human involvement necessary in the initial steps. The user must identify and input the pixel coordinates of matching corners, which becomes prohibitive in models with large numbers of quadrilateral faces. This in turn limits the accuracy of the depth map produced because we are approximating curved objects with piecewise-defined planes. As a large number of small quadrilateral faces yields a better approximation of the original object, the most significant area of future research is in automatically detecting corresponding pixel coordinates.

One possible method is suggested by the projection of the shadow grid onto the hand. The edges and corners of the quadrilaterals tiling the hand are characterized by particularly dark pixels in both photographs. If a pattern of dots were projected rather than a grid, only the corners of the quadrilaterals would be marked by pixels with a distinct intensity. It would then be possible to recover the locations of the pixels with that intensity using \emph{Mathematica}, though steps would have to be taken to ensure the dots were distinct enough from the object to be detected. This would find the locations of the corners in each photograph, but would not necessarily match them. We developed a rudimentary version of this method that used the order of the corner's detection to identify corresponding corners, but this is clearly unreliable, as small rotations of the camera and change in camera angle will change the order in which the corners are detected. An effective matching scheme would be necessary to make this method feasible. Additionally, some way to detect the edges of the object, and to locate the corners of quadrilaterals with one or more edges on the boundary, would be necessary, as the projected dots will not necessarily indicate the edge of the object accurately.

Other areas of future development include making full 3-D models, rather than the bas-relief type models shown here, which would involve combining the depth maps produced by photos from multiple orientations. A related area of research is correcting the distortions from the photograph currently preserved in our models, likely using information presented in \cite{Camera} about length and angle distortion in perspective projections. It is possible that using photographs from multiple orientations will provide enough information to correct the distortions; if not, correcting them will be a necessary element of creating full 3-D models.


\section{Acknowledgments}

This work described in this paper was accomplished during the 2013 REU program at Grand Valley State University.  For more information about 3D printing, including a primer for mathematics professors and their students, please see this web site:  \url{sites.google.com/site/aboufadelreu/}.

This work was partially supported by National Science Foundation under Grant Number DMS-1262342, which funds a Research Experience for Undergraduates program at Grand Valley State University.  Any opinions, findings, and conclusions or recommendations expressed in this material are those of the author(s) and do not necessarily reflect the views of the National Science Foundation.






\bibliography{lite}
\bibliographystyle{mdpi}



%


%

\end{document}